\documentclass[a4paper,12pt]{amsart}

\author{Giles Gardam}
\title{Profinite rigidity in the SnapPea census}

\usepackage[T1]{fontenc}
\usepackage[utf8]{inputenc}

\usepackage{amsfonts}
\usepackage{amsmath}
\usepackage{amssymb}
\usepackage{amsthm}
\usepackage{array}
\usepackage{fullpage}
\usepackage{graphicx}
\usepackage{mathpazo}
\usepackage{microtype}
\usepackage[onehalfspacing]{setspace}
\usepackage{subfig}
\usepackage{textcomp}
\usepackage{tikz-cd}
\usepackage{xcolor}
\usepackage{xspace}

\usepackage{hyperref}

\theoremstyle{plain}
\newtheorem{thm}{Theorem}
\newtheorem{lem}[thm]{Lemma}
\newtheorem{conj}[thm]{Conjecture}
\newtheorem*{strengthenedconj}{Question \ref{qn:profinite_rigidity}'}
\newtheorem*{thmProfinite}{Theorem \ref{thm:my_computation}}
\newtheorem{prop}[thm]{Proposition}
\newtheorem{cor}[thm]{Corollary}
\newtheorem{keythm}{Theorem}

\theoremstyle{definition}
\newtheorem{defn}[thm]{Definition}
\newtheorem{eg}[thm]{Example}
\newtheorem{qn}{Question}

\theoremstyle{remark}
\newtheorem{rmk}[thm]{Remark}

\DeclareMathOperator{\Aut}{Aut}
\DeclareMathOperator{\FC}{FC}
\DeclareMathOperator{\FIA}{FIA}
\DeclareMathOperator{\PGL}{PGL}
\DeclareMathOperator{\PSL}{PSL}

\newcommand{\llim}{\varprojlim}
\newcommand{\normal}{\vartriangleleft}
\newcommand{\Z}{\mathbb{Z}}

\renewcommand{\hat}{\widehat}

\newcommand{\abs}[1]{\lvert #1 \rvert}
\newcommand{\gp}[2]{\langle \, #1 \, | \, #2 \, \rangle}
\newcommand{\setc}[2]{\left\{ #1 \,:\, #2 \right\}}
\newcommand{\Zmod}[1]{{\mathbb{Z}_{#1}}}

\address{Department of Mathematics \\ Technion \\ Haifa \\ Israel}
\email{gilesgar@technion.ac.il}
\thanks{This work was part of the author's doctoral thesis at the University of Oxford.} 

\begin{document}

\begin{abstract}
A well-known question asks whether any two non-isometric finite volume hyperbolic 3-manifolds are distinguished from each other by the finite quotients of their fundamental groups.
At present, this has been proved only when one of the manifolds is a once-punctured torus bundle over the circle.
We give substantial computational evidence in support of a positive answer, by showing that no two manifolds in the SnapPea census of 72 942 finite volume hyperbolic 3-manifolds have the same finite quotients.
\end{abstract}

\maketitle

\section{Introduction}

A standard approach to studying infinite groups is through their finite quotients.
While this has serious limitations in general -- exemplified by the existence of infinite groups having \emph{no} non-trivial finite quotients -- in many contexts, the finite quotients of a group encode much important information about it.
For instance, the fundamental group of any compact 3-manifold is \emph{residually finite}, so it has enough finite quotients that every non-trivial element survives in one.
The question of how much is encoded in the finite quotients of a 3-manifold group has gathered much attention in recent years.
One well-known open question, attributed to Long and Reid in \cite[Question~1]{agol}, is the following:

\begin{qn}
    \label{qn:profinite_rigidity}
    Let $M$ and $N$ be finite volume hyperbolic 3-manifolds.
    If $\pi_1 M$ and $\pi_1 N$ have the same finite quotients, does this imply that $\pi_1 M \cong \pi_1 N$?
\end{qn}
By Mostow Rigidity, $\pi_1 M \cong \pi_1 N$ implies that $M$ and $N$ are isometric.

It is convenient to collect the totality of the finite quotients of a (finitely generated) group $G$ into a single algebraic object, namely its \emph{profinite completion} $\hat{G}$, the inverse limit of the inverse system of its finite quotients, as we introduce in Section~\ref{subsection:profinite}.
This topological group determines the set of (isomorphism classes of) finite quotients, and a well-known result, Proposition~\ref{prop:profinite_from_set}, says that, conversely, it is determined by the set of finite quotients if it is (topologically) finitely generated.

There has been a lot of recent progress in the study of profinite properties of 3-manifolds, with results showing both rigidity and flexibility.
Various properties of 3-manifolds have been shown to be profinite invariants (that is, determined by the profinite completion of the fundamental group), including hyperbolicity, by Wilton--Zalesskii \cite{WZ_geometries}, and being fibred, by Jaikin-Zapirain \cite{jaikin} following \cite{BR_figure_eight}, \cite{BRW_rigidity} and \cite{BF_profinite}.
Following work of Agol and Wise, 3-manifold groups are \emph{good} in the sense of Serre; we content ourselves here with commenting that this implies that the profinite completion detects whether or not $H^3(M)$ is non-trivial, and thus distinguishes cusped from closed hyperbolic 3-manifolds.
The most significant progress on Question~\ref{qn:profinite_rigidity} is a theorem of Bridson, Reid and Wilton \cite{BRW_rigidity}, proving that it holds in the case that $M$ is a once-punctured torus bundle over the circle (so that $\pi_1 M \cong F_2 \rtimes \Z$), building on earlier work that did special cases \cite{BR_figure_eight, BF_profinite}.
The forthcoming paper \cite{eisenstein} gives the first examples of groups which have non-abelian free subgroups and are \emph{absolutely profinitely rigid}, by which we mean that they are each uniquely determined by their profinite completion amongst all finitely generated residually finite groups.
The examples are namely $\PGL(2, \Z[\omega])$ and $\PSL(2, \Z[\omega])$, where $\omega$ is a cube root of unity, which are fundamental groups of hyperbolic 3-orbifolds, both with the figure eight knot complement as a finite sheeted cover.

The purpose of this paper is to report on a computational proof that the manifolds in the benchmark census of (low) finite volume hyperbolic 3-manifolds have pairwise non-isomorphic profinite completions.

\begin{keythm}
    \label{thm:my_computation}
    The 72 942 finite volume hyperbolic 3-manifolds in the SnapPea census are distinguished from each other by the finite quotients of their fundamental groups.
\end{keythm}

These census manifolds are those included in the package \texttt{SnapPy} \cite{SnapPy}, of which 11 031 are closed (available in \texttt{OrientableClosedCensus}) and 61 911 are cusped (available in \texttt{OrientableCuspedCensus}).
The cusped examples represent all orientable cusped hyperbolic manifolds that can be triangulated with at most 9 ideal tetrahedra.

Note that Theorem~\ref{thm:my_computation} does not however answer Question~\ref{qn:profinite_rigidity} in the case that we fix $M$ to lie in the census; our computational method can only prove \emph{relative} profinite rigidity within the census, where \emph{both} $M$ and $N$ must be chosen from the census manifolds.

The value of Theorem~\ref{thm:my_computation} goes beyond the achievement of providing the first wholesale evidence for a positive answer to Question~\ref{qn:profinite_rigidity}.
We mention here one related conjecture, and one direct consequence.

\subsection*{Asymptotic Volume Conjecture}

First, Question~\ref{qn:profinite_rigidity} fits into an important circle of deep work, culminating in the Asymptotic Volume Conjecture following Lück, Bergeron, Venkatesh, Lê, and others.
\begin{conj}
    Let $M$ be a finite volume hyperbolic 3-manifold with fundamental group~$G$.
    Then \[
        \underset{[G:K] \to \infty}{\limsup} \frac{\log \lvert (H_1(K, \Z))_{tor} \rvert}{[G:K]} = \frac{\operatorname{Vol}(M)}{6\pi}.
    \]
\end{conj}
The upper bound on the limit has been proved by Lê \cite{L_bound}.

\begin{figure}
    \centering
    {
    \def\svgwidth{0.9\textwidth}
    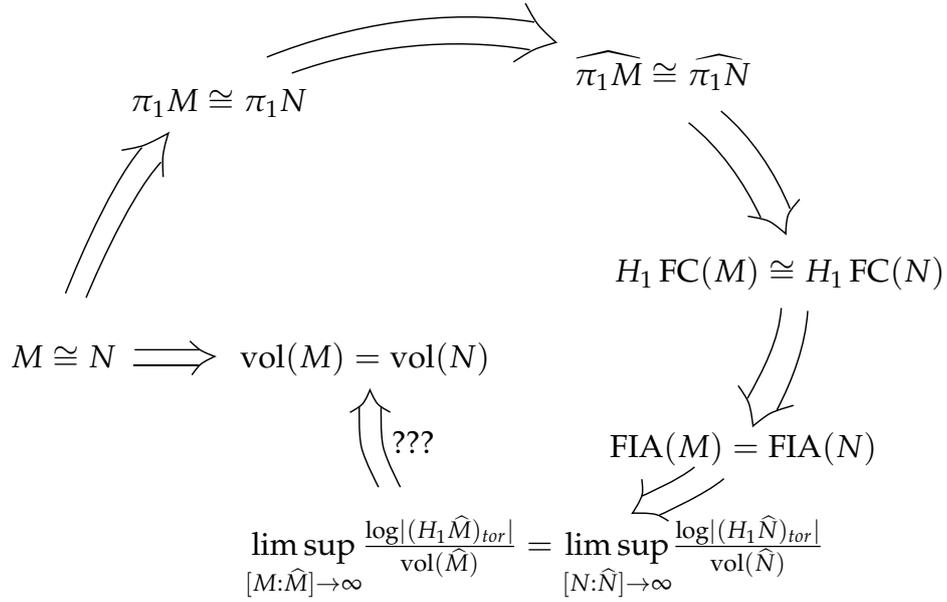
    }
    \caption{Implications for $M$ and $N$ finite volume hyperbolic 3-manifolds}
    \label{fig:spiral}
\end{figure}

One can gradually weaken invariants of manifolds, leading to the spiral of implications shown in Figure~\ref{fig:spiral}.
For a manifold $M$, let $\FC(M)$ denote the lattice of finite sheeted covers of $M$.
We can apply the homology functor (with trivial $\Z$ coefficients) to get $H_1 \FC(M)$, a lattice of abelian groups, which we consider to be annotated by the degrees of the corresponding covers.
A standard fact, which we recount in Section~\ref{subsection:profinite}, says that $\hat{\pi_1 M}$ determines $H_1 \FC(M)$.
If we then forget the lattice information, and the number of subgroups of a given index with a given abelianization, we have \[
    \FIA(M) := \setc{ ([H_1 \check{M}], n) }{ n = [M : \check{M}] < \infty }.
\] That is, $\FIA(M)$ is the set of (isomorphism classes of) abelianizations of finite index subgroups of $\pi_1 M$, together with their indices in $\pi_1 M$.
Many results are consistent with the aphorism that homology of finite sheeted covers of hyperbolic 3-manifolds tell us almost anything we could wish to know about them.
The logical extreme is the following strengthening of Question~\ref{qn:profinite_rigidity}.

\begin{strengthenedconj}
    Let $M$ and $N$ be two finite volume hyperbolic 3-manifolds, and suppose that $\FIA(M) = \FIA(N)$.
    Must $M$ and $N$ be isometric?
\end{strengthenedconj}

It appears that verifying this conjecture in the SnapPea census with current software is infeasible; see Remark~\ref{rmk:fia_failure} below.

The implications in the figure are drawn as a spiral because, conjecturally, we end up back where we started up to ``finite ambiguity'', since only finitely many hyperbolic 3-manifolds can have the same finite volume \cite[3.6~Theorem]{thurston_bams}.

We note that $\FIA(M)$ is a strictly weaker invariant than $\hat{\pi_1(M)}$, even when $M$ is a compact 3-manifold (although conjecturally this cannot happen in the hyperbolic case).
The following example is due to Gareth Wilkes.
\begin{eg}
    \label{eg:wilkes}
    There are Seifert fibred 3-manifolds $M$ and $N$ such that $\FIA(M) = \FIA(N)$ but $\hat{\pi_1 M} \not \cong \hat{\pi_1 N}$.
    Indeed, the fundamental groups
    \begin{align*}
        \pi_1 M &= \gp{a, b, c, h}{h \text{ central},\, a^4 h,\, b^4 h,\, c^2 h,\, abc} \\
        \pi_1 N &= \gp{a, b, c, h}{h \text{ central},\, a^4 h^3,\, b^4 h^3,\, c^2 h,\, abch}
    \end{align*}
    are distinguished by their maximal 2-class 6 quotients (of order $2^{12}$).
    These manifolds are commensurable.
\end{eg}
As these are Seifert fibred 3-manifolds (and not of the form of Hempel's counterexamples), we can apply \cite[Theorem~1.2]{gareth} to conclude \emph{a priori} that $\pi_1 M \not \cong \pi_1 N$ implies that $\hat{\pi_1 M} \not \cong \hat{\pi_1 N}$.

\subsection*{No duplicates in the census}

Second, Theorem~\ref{thm:my_computation} gives independent verification that the census does not contain any duplicates.
The standard way of verifying this is to compute the canonical Epstein--Penner cell decomposition.
However, rounding errors in imprecise computational arithmetic of real numbers has previously lead to duplicates.
One such pair was identified by Burton \cite{B_duplicate}.
Our verification -- while of course dependent on many large computer calculations that cannot be replicated by hand -- involves only precise discrete computations, in combinatorial group theory and linear algebra over $\Z$.

Nathan Dunfield has informed us that he is applying similar techniques to distinguish closed hyperbolic 3-manifolds for a new, extended census.

The outline of the paper is as follows. In Section~\ref{section:theory} we introduce some basic theory of profinite completions and hyperbolic 3-manifolds, then in Section~\ref{section:practice} we discuss practical matters: heuristics, limitations, and our methods. Section~\ref{section:results} gives the results, and future work is discussed in Section~\ref{section:future}.

\section{Theory}
\label{section:theory}

If we proceed on the assumption that the census manifold \emph{are} distinguished from each other by the finite quotients of their fundamental groups, there is a naive algorithm that will prove this for us: simply enumerate finite groups, and determine which manifold groups have them as quotients, until we have enough information to distinguish them all.
However, in general we cannot give any prediction as to how long such a verification would take, as discussed below.
Indeed, there is a common expectation that the time needed to prove Theorem~\ref{thm:my_computation} would be astronomical (for instance, the 150 groups with trivial abelianization were not distinguished from each other by counting maps onto finite simple groups after several weeks of computation).
However, structure theory of the profinite completion and its subgroups and some theory of hyperbolic 3-manifolds reveal why a less naive approach (not simply enumerating maps to finite groups) should be feasible, as our computations have demonstrated.

\subsection{Computing with the profinite completion}

\label{subsection:profinite}

We refer to \cite{RZ} for fundamental results on profinite completions.

\begin{defn}
    Let $G$ be a group.
    The \emph{profinite completion} $\hat{G}$ of $G$ is the inverse limit \[
        \llim_{N \normal G, [G:N] < \infty} G / N
    \] of the inverse system of finite quotients of $G$.
\end{defn}

The profinite completion has the expected universal property.
\begin{lem}
    There is a natural map $\iota \colon G \to \hat{G}$ such that every map from $G$ to a finite group $Q$ factors through $\hat{G}$ uniquely.
    \[
    \begin{tikzcd}
    {} & \hat{G} \arrow[dashed]{d}{\exists ! \Phi} \\
        G \arrow[hook]{ru}{\iota} \arrow{r}{\phi} & Q
    \end{tikzcd}
    \]
\end{lem}

The map $\iota$ is an embedding if and only if $G$ is residually finite.

\begin{prop}[{\cite[Corollary~3.2.8]{RZ}}]
    \label{prop:profinite_from_set}
    Let $G_1$ and $G_2$ be finitely generated groups with the same finite quotients.
    Then $\hat{G_1} \cong \hat{G_2}$.
\end{prop}

\begin{lem}[{\cite[Proposition~3.2.2]{RZ}}]
    Let $H$ be a finite index subgroup of $G$.
    Let $\overline{H}$ denote the closure of $\iota(H)$ in the profinite topology on $\hat{G}$.
    Then $\overline{H} \cong \hat{H}$, and the isomorphism is natural.
\end{lem}

\begin{cor}
    There is a one-to-one correspondence between finite index subgroups of a finitely generated group $G$ and finite index (open) subgroups of $\hat{G}$.
    This bijection preserves profinite completion, and thus abelianization, as well as normality and the isomorphism class of quotients by normal subgroups.
\end{cor}
The word \emph{open} can be omitted in the statement of the corollary, by the Nikolov--Segal Theorem \cite{NS_profinite} (for which finite generation of $G$ is essential); this is not needed for our applications of this fact, since we only compare finite index subgroups of discrete groups whose profinite completions are isomorphic.

By combining the universal property, and Proposition~\ref{prop:profinite_from_set}, we see that the set of (isomorphism classes of) finite quotients of a finitely generated group $G$ determines in particular the \emph{number} of surjections of $G$ onto any finite group.

Two standard approaches to proving a group $G$ to be non-trivial are
\begin{itemize}
    \item find a non-trivial finite group onto which $G$ maps (or equivalently, a proper finite index subgroup); \emph{and}
    \item show that $G$ has non-trivial abelianization.
\end{itemize}
To show that a group is infinite, one combines these approaches and attempts to find a finite index subgroup of infinite abelianization.

These techniques are actually two sides of the same coin: both are attempting to compute invariants of the profinite completion $\hat{G}$.

In contrast to the majority of decision problems for groups that one encounters, for profinite completions of finitely presented groups, distinguishing groups $\hat{G}$ and $\hat{H}$ is the easy direction: one can enumerate all maps from the groups to finite groups (for instance, via maps to finite symmetric groups of increasing degree) and wait until there is a finite group onto which $G$ maps, say, but $H$ doesn't.
The other direction is unsolvable: Bridson and Wilton proved that one cannot even decide if the profinite completion of a finitely presented group is trivial \cite{BW_profinite_triviality}.
Thus there can be no general computable bound on the time needed for an algorithm to provide proof that two groups have different profinite completions when that is the case, because this would allow us to conclude that the profinite completions are isomorphic, if this is the case, after the allowed computational time has elapsed.

Suppose that one has a set $S$ of groups which is \emph{relatively profinitely rigid}; that is, for $G, H \in S$, we have $G \cong H$ if and only if $\hat{G} \cong \hat{H}$.
Then given two finite presentations of groups in $S$ we can decide whether the two groups are isomorphic: by day, we attempt to construct a proof that $G \cong H$, and by night, we attempt to prove that $\hat{G} \not \cong \hat{H}$.
It is guaranteed that one of these procedures will terminate:
\begin{prop}
    \label{prop:isom_prof_rigid}
    The isomorphism problem is solvable in a relatively profinitely rigid class of finitely presented groups.
\end{prop}

Note that the profinite completion does \emph{not} determine the nilpotent quotients of a group \cite{remeslennikov_nilpotent}.
Thus, a commonly used computational technique for distinguishing groups is not available to us.

\subsection{Hyperbolic 3-manifolds}

Non-abelian simple quotients have previously proved to be effective at distinguishing profinite completions of groups.
This was a natural place to look in work on parafree groups of Baumslag--Cleary--Havas \cite{BCH_exp}, because nilpotent quotients cannot distinguish those groups.
For the hyperbolic 3-manifolds in the SnapPea census, nilpotent quotients are often similarly ruled out because $b_1(M, \Z/p) \leq 1$ for all primes (for example, for many one-cusped examples, such as all knot complements, with $H_1(M, \Z) \cong \Z$), which means no non-abelian nilpotent quotients.
The profinite completion already encodes all the abelian quotients.

A further theoretical justification for using simple groups is that Long--Reid proved that finite volume hyperbolic 3-manifolds are residually simple \cite{LR_simple}.

Dunfield and Thurston tabulated for each finite simple group of order up to 32~736 the percentage of closed census manifolds with that simple group as a quotient of the fundamental group \cite[Table~2, p.~12]{DT_haken}.
Amongst these finite simple groups and these manifolds, a random manifold has a random simple group as a quotient with probability 34.7\%.
The Mathieu group $M_{11}$ is a quotient of only 17.1\% of the manifolds, which is minimal for the simple groups considered.

Note that the achievement of Bridson--Reid--Wilton was to show that the set of groups $F_2 \rtimes \Z$ is profinitely rigid.
It is \emph{not} known whether extensions of non-abelian free groups by cyclic groups are profinitely rigid in general.
Many examples, in fact more than 47\% of the census manifolds, are free-by-cyclic (as determined by Brown's Criterion \cite[\S~4]{brown}).
The resolution of the Virtual Fibering Conjecture by Agol, following work of Wise and coauthors, means that any finite volume cusped hyperbolic 3-manifold is virtually free-by-cyclic.

A more ambitious variant of Question~\ref{qn:profinite_rigidity} would be to allow one of the two groups to be any finitely generated residually finite group, that is, to ask whether finite volume hyperbolic 3-manifolds are absolutely profinitely rigid; it is open whether this greatly strengthened conjecture holds.
Remesslennikov's question of whether free groups are distinguished in the class of all finitely generated residually finite groups by their profinite completions also remains open.

\section{Practice}
\label{section:practice}

\subsection{Heuristics}

Two strategies to distinguish profinite completions are
\begin{itemize}
    \item find abelianizations of finite index subgroups; \emph{and}
    \item count maps onto finite (simple) groups (up to automorphisms of the quotient).
\end{itemize}
By factoring out automorphisms of the quotient, we are counting the normal subgroups which give the specified finite group as quotient.

One would hope to distinguish hyperbolic 3-manifolds using $\FIA$.
Certainly, it appears experimentally that this approach is more effective that enumerating finite quotients.
This mimics the phenomenon of Dunfield--Thurston \cite{DT_haken}, where it was easier to verify the virtual Haken conjecture by verifying a stronger, more algebraic, result, namely virtual positive first Betti number, which also involves abelianizations of finite index subgroups.
A fundamental difference in character between the virtual Haken conjecture and Question~\ref{qn:profinite_rigidity} is that one only needs to exhibit a Haken cover to prove the conjecture in a specific instance.
In our present case, not only are we unable to prove profinite rigidity relative to the class of all finite volume hyperbolic 3-manifold groups, it is difficult to imagine an easily verifiable certificate just within the census, because distinguishing profinite completions requires proving \emph{non-existence} of certain quotients or certain abelianized subgroups, and it is not at all clear how one might do this without repeating the exhaustive enumeration.

However, enumerating all subgroups of index $n$ has complexity factorial in $n$, and many manifolds have the same abelianizations of low index subgroups.

We thus turn to a powerful combination of the two naive strategies: determining maps onto finite simple groups, then computing the abelianization of the kernel.
There is a very good heuristic reason for why this should be effective; we quantify this using entropy, the expected information in a random variable.

\begin{defn}
    Let $X$ be a discrete random variable taking values in $\{x_1, \dots, x_n\}$, each with probability $P(x_i)$.
    The \emph{entropy} of $X$ is \[
        H(X) := - \sum_{i=1}^n P(x_i) \log_2 P(x_i).
    \]
\end{defn}

The maximal entropy of a random variable taking $n$ values is $\log_2 n$ bits, and occurs when it is the uniform distribution (each value occurs with probability~$\frac{1}{n}$).
Our random variables will be the value(s) of invariants for the $N = 72942$ census manifolds, drawn uniformly.
The entropy is $\log_2 N$ if and only if a given set of invariants completely distinguishes the manifolds.

The number of maps from a group $G$ defined by a $k$ generator presentation onto a finite group $Q$ of order $n$ is certainly at most $n^k$.
(For non-abelian simple groups, with few outer automorphisms and very high probability that a random pair generates, $n^{k-1}$ is a fairly good approximation for the upper bound coming from the free group.)
Thus, computing the number of maps to a given finite group can only provide us with approximately $k \log_2 n$ bits of information.

\begin{figure}
    \centering
    \subfloat{{\includegraphics[width=0.45\textwidth]{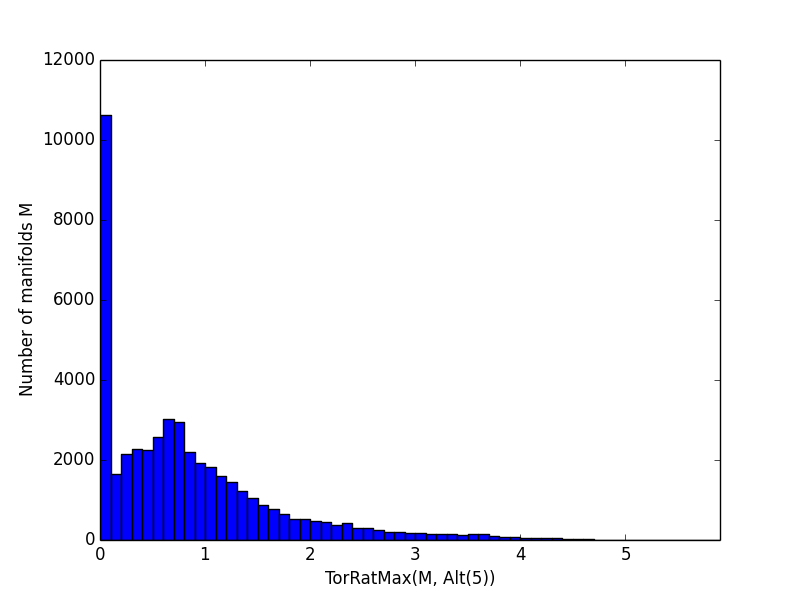}} }%
    \quad
    \subfloat{{\includegraphics[width=0.45\textwidth]{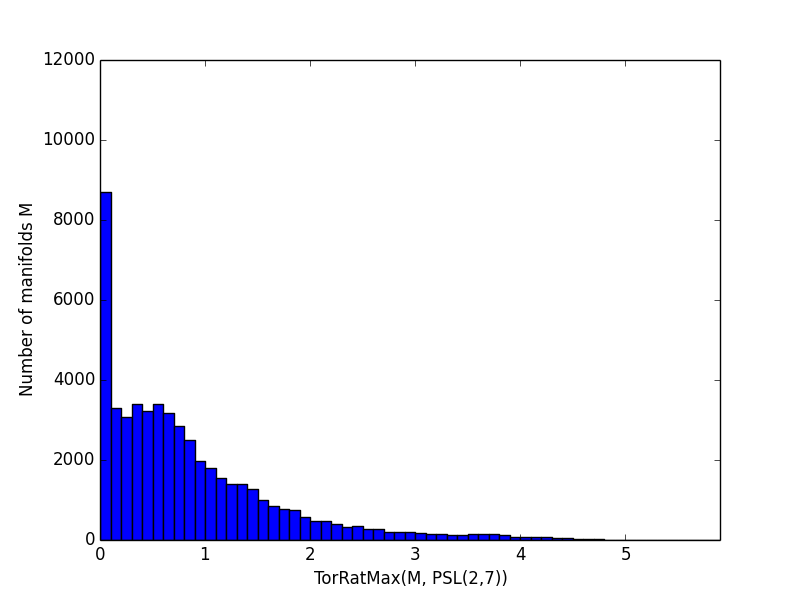}} }%
    \quad
    \subfloat{{\includegraphics[width=0.45\textwidth]{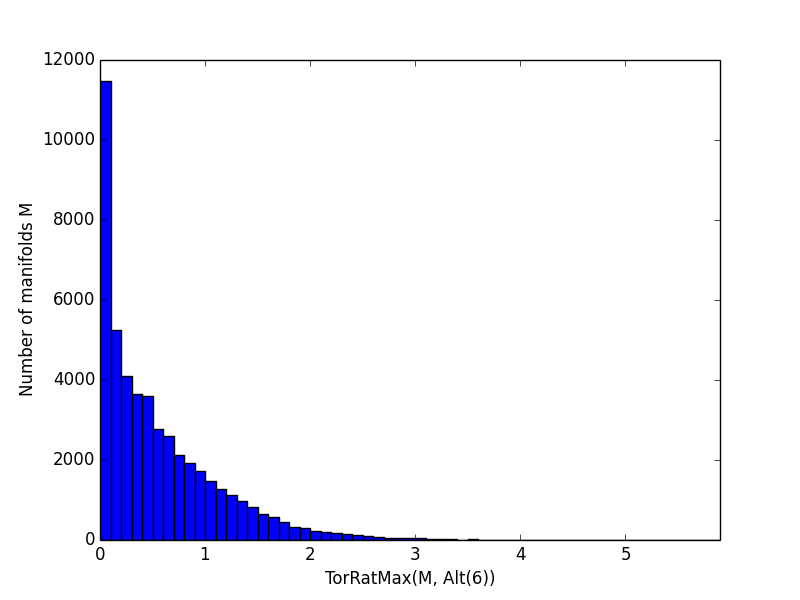}} }%
    \caption{Histogram for $\max \setc{ \frac{6 \pi \log \lvert (H_1(K, \Z))_{tor} \rvert}{\lvert Q \rvert \operatorname{Vol}(M)}}{ K \text{ s.t. } \pi_1 M / K \cong Q}$ amongst manifolds $M$ having a cover with Galois group $Q$, which accounts for at least 64\% of the census manifolds in each of the cases $Q \cong A_5, \PSL(2,7), A_6$ illustrated. The Asymptotic Volume Conjecture predicts a value of 1.0 for the limit supremum as the index goes to infinity.}
    \label{fig:AVC} 
\end{figure}

On the other hand, torsion in homology grows very quickly (the torsion-free rank also provides useful, albeit secondary, information).
We expect its size, following the Asymptotic Volume Conjecture, to be on the order of $e^{n \operatorname{Vol}(M) / 6 \pi}$ for a normal subgroup of index $n$.
Figure~\ref{fig:AVC} on page~\pageref{fig:AVC} gives an illustration of how closely the actual homology groups computed for certain simple covers correspond with the predictions of the conjecture.
Thus we get entropy that is linear, rather than logarithmic, in the size of the quotient $Q$ considered, provided the homology groups arising are sufficiently varied.
(In considering these heuristic estimates, we must keep in mind that there is a bound of $\log_2 N$ on entropy for any collection of invariants, and that the random variables will not be independent so we cannot simply sum up their respective entropies.)
For any given $Q$, a large but not overwhelming percentage of manifolds will have the same amount of information under both schemes, as they have no surjections onto $Q$; when a group has more than 1 surjection to $Q$, this only helps us.

\begin{rmk}
There are various ways one could measure how a set of invariants contributes towards distinguishing a finite set of objects.
We believe that entropy is the best measure.
Simple alternatives, such as counting the number of equivalence classes or the number of objects that have been distinguished from all the others, fail to capture the ``shape'' of the partition.
Consider for example two possible partitions of 10 objects, either as $5$ pairs or as sextuple and $4$ singletons.
Entropy captures well the intuitive view that partitioning the set into $5$ pairs is better progress.
\end{rmk}

Table~\ref{table:entropies} lists the entropy of counting maps to some finite simple groups for the sample space of SnapPea census manifolds (taken uniformly at random).
The last two columns demonstrate that a lot of the gap between entropy of the homology of covers $K$ and the upper bound of $\log_2 (72 942)$ is accounted for by the number of manifolds which have no such cover.

\begin{table}
\centering
\begin{tabular}{c | r r c c}
\multicolumn{1}{c}{} & \multicolumn{2}{c}{entropy of:} & \multicolumn{2}{c}{when $\#\{K\} \geq 1$:}\\
\multicolumn{1}{c}{group} & \multicolumn{1}{c}{$\#\{K\}$} & \multicolumn{1}{c}{$\{H_1 K\}$} & entropy & $\log_2(\# \text{manifolds})$\\
\hline
  $A_5$     & 2.37 &  8.90 & 12.37 & 15.51 \\
$\PSL(2, 7)$ & 2.88 & 10.82 & 13.98 & 15.66 \\
 $A_6$      & 2.80 &  9.92 & 13.53 & 15.57
\end{tabular}

\caption[Entropy in maps to simple groups]{Entropy in the number of regular covers $K$ of $M$ with Galois group isomorphic to a given finite simple group, for $M$ in the census; entropy in the (multi)set of abelianizations of such covers; entropy amongst only those manifolds with at least one such simple cover; upper bound on that entropy}
\label{table:entropies}
\end{table}

We now recount a very concrete example of the power of computing abelianizations of kernels.
The two knots identified in SnapPy as \texttt{K14a3482} and \texttt{K14a3494} are very difficult to distinguish.
For full reference, their Dowker--Thistlethwaite codes are

{
\begin{center}
\texttt{4 10 14 16 2 24 22 18 8 6 26 28 20 12} \quad \text{and} \\
\texttt{4 10 14 16 2 26 24 18 8 6 12 28 20 22} \quad \phantom{\text{and}}
\end{center}
}

\noindent respectively.
Their complements have the same volume 24.1942..., their Alexander polynomials are identical, they have the same Khovanov homology (computed with \texttt{KnotKit} \cite{knotkit}), and they have the same number of surjections onto any simple group of order less than 2500.
Nonetheless, the abelianization of kernels of maps to the simple group $\operatorname{PSL}(2, 7)$ suffices to show that their profinite completions differ.

As well as having high entropy, the kernels (as a set) are characteristic, and the problem of proving non-existence is dissolved into enumerating maps to a finite quotient.
Enumerating non-normal subgroups even of index $60 = \lvert A_5 \rvert$ is completely infeasible.

\subsection{Difficulties and limitations}

The main concerns at this point are twofold: the gap between predictions of torsion in homology and actual low volume results, and the correlation between different invariants.

A practical difficulty, which we are yet to explain satisfactorily, is that GAP has extreme difficulty computing \texttt{GQuotients} on a small number of examples (less than 1 in 1000, at least for the smallest of simple quotients).
For instance, the fundamental group $G$ of the manifold identified in SnapPy as \texttt{t05599(0,0)} has the presentation
\[
    \gp{a,b,c}{a^2 b^5 a^2 b^2 c^{-2},\, a^5 c^3 b^{-2}}.
\] It has one normal subgroup $N$ of index 60 with quotient $G / N \cong A_5$, and $H_1 N \cong \Z^{12} \times \Z/12$.
Computing all the surjections from $G$ to $A_5$ with \texttt{GQuotients} (which works very directly with the finite presentation) takes GAP 20 minutes and requires gigabytes of memory.
We rolled our own method, enumerating all maps from the group $G$ to its quotient $Q$.
The number of iterations we run is bounded by $\abs{\operatorname{Inn}(Q) \backslash Q} \cdot \abs{Q}^{d(G)-1}$, because we can assume without loss of generality that the first generator of $G$ is sent to a preferred element in each conjugacy class, and then send the second generator to a preferred element modulo the centralizer of the image of the first element, and so on.
In short, the tuple of images of generators is chosen to be minimal under the lexicographic ordering (after picking an arbitrary order on $Q$) within its conjugacy class.
We simply check at the end which maps give the same kernel; since simple groups have such small outer automorphism groups, we will by this point have only overcounted by a factor of 2 or 4 usually (and it is not worth the hassle of explicitly factoring out by the action of $\Aut(Q)$ and not just $\operatorname{Inn}(Q)$ before this stage).
Our method computes all surjections of the aforementioned $G$ onto $A_5$ in under a second, a speed up of over 1000.

\subsection{Methods}

We computed the minimum number of invariants to distinguish the groups.
That is, as soon as a group had been distinguished from the others, it was removed from consideration.
We computed
\begin{itemize}
    \item abelianizations of subgroups up to index 7;
    \item the abelianization of the maximal abelian cover if that was finite index, and failing that the cyclic covers up to index 10 if the abelianization was $\Z$; \emph{then}
    \item abelianizations of kernels of maps to small finite non-abelian simple groups
\end{itemize}
This was performed in parallel in twenty cores, coordinated by a \texttt{python} script.
We used the wonderful program \texttt{SnapPy} \cite{SnapPy} to work with the manifolds in question, and in particular to extract group presentations, and used \texttt{GAP} \cite{GAP4} for all the group theoretic computations.

\section{Results}
\label{section:results}

\begin{thmProfinite}
    The 72 942 finite volume hyperbolic 3-manifolds in the SnapPea census are distinguished from each other by the finite quotients of their fundamental groups.
\end{thmProfinite}

This took around 64 hours of CPU time.

\begin{figure}
    \centering
    \subfloat[Number of manifolds distinguished]{{\includegraphics[width=0.45\textwidth]{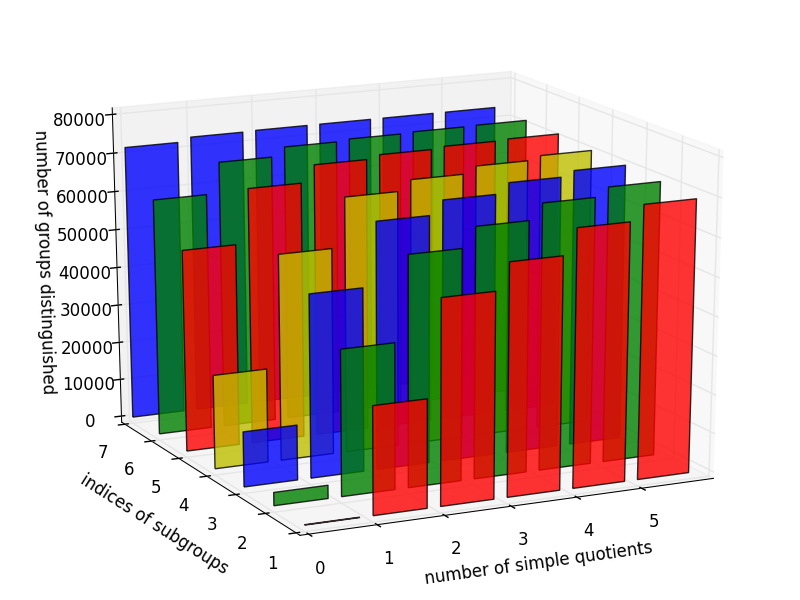} }}%
    \quad
    \subfloat[Entropy]{{\includegraphics[width=0.45\textwidth]{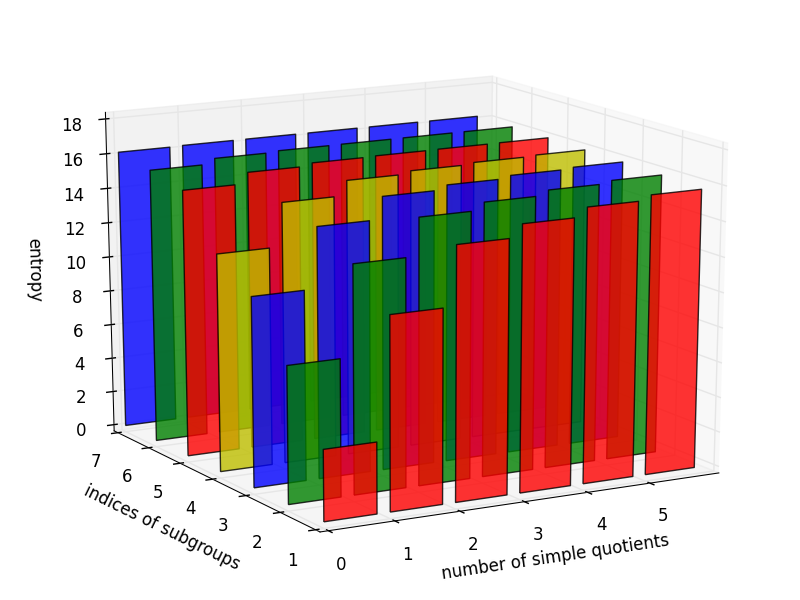} }}%
    \caption{Information in low-index subgroups and small simple quotients. Each bar corresponds to the information contained in both the abelianization of all subgroups up to a given index and the abelianization of kernels of maps onto the smallest 1, 2, 3, 4 or 5 non-abelian simple groups (hence the monotonicity).}
    \label{fig:entropy_etc}
\end{figure}

A plot indicating the number of manifolds distinguished and the entropy from computing abelianization of
\begin{itemize}
    \item all subgroups up to a given index, together with
    \item all kernels of maps onto the smallest 1, 2, 3, 4 or 5 non-abelian simple groups
\end{itemize}
is indicated in Figure~\ref{fig:entropy_etc}.
The plots are very similar; we note however that the integral homology of a manifold alone, which only distinguishes 102 or the manifolds, still has approximately 4 bits of information.

\begin{rmk}
\label{rmk:fia_failure}
We were unsuccessful in distinguishing the groups using only $\FIA$.
For the 11 031 closed manifolds in the census, two months of CPU time was insufficient.
At this point, there were 13 manifolds left to be differentiated, in 5 pairs and 1 triple.
No proper subgroups had been found for these 13 examples up to index 12.
Beyond this index, the exhaustive search for finite index subgroups can take months for a single group.
Under the reasonable assumption that these groups do have enough subgroups waiting at index 13, verifying that $\FIA$ distinguishes them would take years of CPU time.
Thus, it is unlikely that we could succeed at this task in our lifetime without implementing a parallel algorithm for enumerating low-index subgroups.
\end{rmk}

A very natural question to ask at this point, especially in light of the preceding remark, is: what indices of subgroup are needed, and what order of simple quotients? For instance, a group with no subgroup of index up to 12 cannot map onto any of the 5 simple groups of order less than 1000, because they all have low index subgroups that could be pulled back.

The answer is that the largest simple quotient used was $\PSL(2, 23)$, of order 6072.
We tabulate the number of manifolds groups whose profinite completions had been distinguished from all the others at each stage of the computation in Table~\ref{table:simples}.
In addition, for the 13 non-abelian simple groups used we list their order and minimal index of a proper subgroup.
With the exception of $\operatorname{PSU}(3, 3)$ (which happened to provide no useful information), they all have a subgroup of index 24 or less.

\begin{rmk}
    This does not mean that all of the 3-manifold groups considered have a subgroup of index 24 or less: we stopped computing subgroups of a group as soon as it was distinguished from all the others.
    What we can definitely say is that every group has a subgroup of index at most 401, the largest prime $p$ such that there exists a manifold $M$ for which $\Zmod{p}$ is the smallest non-trivial quotient of $H_1 M$.
    The only such $M$ for $p = 401$ is \texttt{v1860(2,3)} (with $H_1 M \cong \Zmod{401}$), so we did not need to find any other subgroups of $\pi_1 M$ (it does, however, surject onto $\PSL(2, 14)$, so it has a subgroup of index 14).
    The last of the 150 perfects groups remaining was the fundamental group of \texttt{s636(-4,3)}, which was distinguished by virtue of having no maps onto the non-abelian simple groups up to $\PSL(2, 16)$.
    We determined separately that it maps onto $\PSL(3, 3)$, so it has a subgroup of index 13 (much less than 401).
\end{rmk}

\begin{table}
\centering
\begin{tabular}{r r}
\multicolumn{1}{c}{Invariant} & \multicolumn{1}{c}{\# dist.} \\
\hline
FIA to index 1 &    102  \\
FIA to index 2 &   3317  \\
FIA to index 3 &  10837  \\
FIA to index 4 &  10095  \\
FIA to index 5 &  28068  \\
FIA to index 6 &   9217  \\
FIA to index 7 &  10029  \\
abelian covers &    966  \\
\end{tabular}
\qquad
\begin{tabular}{r r r r}
\multicolumn{1}{c}{Group} & \multicolumn{1}{c}{Order} & \multicolumn{1}{c}{Min. index} & \multicolumn{1}{c}{\# dist.} \\
\hline
$A_5$       & 60   & 5  &   8 \\
$\PSL(2, 7)$  & 168  & 7  &  12 \\
$A_6$       & 360  & 6  &   4 \\
$\PSL(2, 8)$  & 504  & 9  & 101 \\
$\PSL(2, 11)$ & 660  & 11 &  82 \\
$\PSL(2, 13)$ & 1092 & 14 &  51 \\
$\PSL(2, 17)$ & 2448 & 18 &  37 \\
$A_7$       & 2520 & 7  &   0 \\
$\PSL(2, 19)$ & 3420 & 20 &   8 \\
$\PSL(2, 16)$ & 4080 & 17 &   2 \\
$\PSL(3, 3)$  & 5616 & 13 &   2 \\
$\operatorname{PSU}(3, 3)$  & 6048 & 28 &   0 \\
$\PSL(2, 23)$ & 6072 & 24 &   4 \\
\end{tabular}

\caption[Invariants computed and simple groups]{The number of groups distinguished at each stage, and for the 13 smallest non-abelian simple groups also the order and minimal index of proper subgroups}
\label{table:simples}
\end{table}

\section{Future work}
\label{section:future}

\subsection*{Technical improvements}

If one were to attempt to distinguish the manifolds using $\FIA$ alone, it would be possible to construct a partial certificate, from which one could reliably re-prove the appropriate version of Theorem~\ref{thm:my_computation} more quickly.
For instance, if at some point it was necessary to distinguish a group $G$ with a $2$-generator presentation and a group $H$ with a $3$-generator presentation, a certificate could record a description of a subgroup of $H$ (as stabilizer of a point in a permutation representation), and then one only needs to show that $G$ has no subgroup at that index with that abelianization; this is much faster than enumerating subgroups in $H$ exhaustively, because of the difference in presentation rank.
We emphasize that the time spent finding low index subgroups dominates the time spent computing their abelianizations.

Many invariants turned out to provide little entropy (or not to help distinguish some subset of the groups at all).
An interesting and natural question is: given all the invariants computed, and the times of computation (or re-computation, modulo a potential certificate) per invariant per group, what is a collection of pairs of invariant and group, whose determination proves that the groups are distinct, with minimal computation time?
This is in fact an NP-complete problem, as was proved by Hyafil and Rivest \cite{HR_binary}.
However, given that, especially for low index subgroup enumeration, the time required grows quickly, it seems that heuristics would allow for a very good approximation to the optimal binary decision tree.

One could also exploit the lattice structure of $\FIA$, although thus far we have not produced a concrete example where this provides additional information over the multiset of abelianizations together with subgroup index.
The extra complication of solving the isomorphism problem for lattices is a dissuading factor, at least in implementing the use of this additional structure (we imagine that the time spent computing the subgroups would however still dominate the time solving the lattice isomorphism problem).

\subsection*{Knot groups}

One is most interested in profinite rigidity in settings where the fundamental group determines the manifold.
This is not true of knot complements in general, but it is true for complements of prime knots (modulo mirror symmetry).
Thus Boileau and Friedl proposed \cite{BF_profinite} the question of whether the complements of prime knots are profinitely rigid.
From an experimental point of view, however, this question is almost exactly the same as Question~\ref{qn:profinite_rigidity}: of the 1 701 935 prime knots of crossing number at most 16 tabulated by Hoste, Thistlethwaite and Weeks \cite{onepointseven}, only 32 are non-hyperbolic.
We plan to apply our machinery to this collection of examples in any case.
Fortunately, the knots are available in SnapPy with group presentations of quite low rank: 67\% are rank 3 and 31\% are rank 4, with only 10 examples of rank 6, the largest occurring.
(For comparison, the standard Wirtinger presentation for a knot group given a knot diagram has as many generators as the diagram has crossings, which is prohibitive when equal to 16.)

Moreover, we are lucky that the unsigned Alexander polynomial (that is, the Alexander polynomial modulo multiplication by $\pm 1$) has recently been shown to be a profinite invariant \cite{ueki}.
This takes us most of the way: we computed the entropy of the unsigned Alexander polynomial amongst the prime knots of crossing number at most $16$ to be $16.67$, which is very close to the entropy of $\log_2 (1701935) = 20.70$ that full discrimination of these knots would require.
An interpretation of the gap of $4.03$ is that a random knot lies in an unsigned Alexander polynomial equivalence class of size $2^{4.03} \approx 16$, where this average is computed as the geometric mean.
Indeed, 49.6\% of the knots are in an equivalence class of size at most 16.
There are 140 261 knots, representing 8.2\% of the total, that are each already distinguished from all the others by their unsigned Alexander polynomial, which leaves us with 1 561 674 knots.

Thus we are optimistic that by exploiting some knot theory we will be able to carry out this experiment to completion.
It will, however, take much longer than the SnapPea census experiment.
In particular, a random sample of 1000 of the 57 005 pairs of knots with the same Alexander polynomial (modulo sign) took approximately 4.5 hours of CPU time.
At this rate, extrapolating generously -- whereas we anticipate that the larger Alexander polynomial equivalence classes would take longer \emph{per knot} -- we have an estimated lower bound of 150 days of CPU time.
By a theorem of Fox, the torsion in homology of the cyclic covers of a knot is determined by its Alexander polynomial.
Moreover, the homology of cyclic covers will be periodic when all roots of the Alexander polynomial are roots of unity \cite{gordon}, so we cannot get away with using just these obvious subgroups, and will inevitably have to search for non-nilpotent finite quotients, and most likely finite simple quotients.
Each such \texttt{GQuotients} search takes on the order of seconds; 1.5 million seconds is approximately 17 days.

\subsection*{Acknowledgements}
I thank Martin Bridson for guidance and support while supervising this project as part of my thesis, Alan Reid for a stimulating minicourse at YGGT~VI in Oxford (March 2017) which lead to me undertaking this project, and Nathan Dunfield, Jim Howie, Marc Lackenby, Nikolay Nikolov, and Gareth Wilkes for their comments.

\bibliographystyle{alpha}
\bibliography{profinite}

\end{document}